\documentclass{article}

\usepackage{amsmath}
\usepackage{amssymb}
\usepackage{latexsym}
\usepackage{amsxtra}
\usepackage{amscd}
\usepackage{theorem}
\usepackage{xy}
\usepackage[utf8]{inputenc}
\usepackage{hyperref}
\usepackage[usenames,dvipsnames]{color} 
\usepackage{tikz}
\numberwithin{equation}{section}
\numberwithin{table}{section}
\numberwithin{figure}{section}

\usepackage{booktabs} % michel, 25.3.
\theoremstyle{plain}

{\theorembodyfont{\slshape}

        \newtheorem{thm}{Theorem}[section]
        \newtheorem{cor}[thm]{Corollary}
        \newtheorem{lem}[thm]{Lemma}
        \newtheorem{prop}[thm]{Proposition}

}

{\theorembodyfont{\rmfamily}

        \newtheorem{rem}[thm]{Remark}
        
        \newtheorem{exa}[thm]{Example}
        
        \newtheorem{notation}[thm]{Notation}

}

\renewcommand{\em}{\sl}

\newcommand{\proof}{{\bf Proof:\ }}
\newcommand{\Endproof}{\hspace{\fill} $\Box$\vspace{1ex}\noindent}

\makeatletter
\renewcommand{\subsection}{\@startsection{subsection}{2}%
        {\z@}{-3.25ex plus -1ex minus-.2ex}{-1em}{\bf}}
\makeatother

\newcommand{\NN}{\mathbb{N}}
\newcommand{\ZZ}{\mathbb{Z}}
\newcommand{\QQ}{\mathbb{Q}}
\newcommand{\FF}{\mathbb{F}}

\newcommand{\PP}{\mathbb{P}}

\newcommand{\Gal}{\mathop{\rm Gal}\nolimits}

\newcommand{\p}{\mathfrak{p}}
\newcommand{\OO}{\mathcal{O}}
\newcommand{\Spec}{\mathop{\rm Spec}}

\newcommand{\Aut}{\mathop{\rm Aut}\nolimits}

\newcommand{\Frob}{\mathop{\rm Frob}\nolimits} 
\newcommand{\Ind}{\mathop{\rm Ind}\nolimits}

\newcommand{\gen}[1]{\mathopen\langle#1\mathclose\rangle}
\newcommand{\abs}[1]{\lvert#1\rvert}

\newcommand{\X}{\mathcal X}
\newcommand{\Y}{\mathcal Y}

\newcommand{\Sc}{\mathcal S}
\newcommand{\Xb}{\bar{X}}
\newcommand{\Yb}{\bar{Y}}
\newcommand{\Zb}{\bar{Z}}
\newcommand{\xb}{\bar{x}}
\newcommand{\yb}{\bar{y}}
\newcommand{\ub}{\bar{u}}
\newcommand{\vb}{\bar{v}}

\newcommand{\et}{{\rm et}}
\newcommand{\gb}{\bar{g}}
\newcommand{\hb}{\bar{h}}
\newcommand{\fb}{\bar{f}}

\begin{document}

\title{The functional equation for $L$-functions of  hyperelliptic curves}

\author{Michel B\"orner, Irene I.\ Bouw, and Stefan Wewers}
\date{}
\maketitle

\begin{abstract}
  We compute the $L$-functions of a large class of algebraic curves, and
  verify the expected functional equation numerically. Our computations are
  based on our previous results on stable reduction to calculate the local
  $L$-factor and the conductor exponent at the primes of bad reduction.  Most
  of our examples are hyperelliptic curves of genus $g\geq 2$ defined over
  $\QQ$ which have semistable reduction at every prime $p$.  We also treat a
  few more general examples of superelliptic curves.

\vspace{1ex}\noindent 
2010 {\em Mathematics Subject Classification}. Primary 11G40. 
Secondary: 14G10, 11G20.
\end{abstract}

\section{Introduction}

\subsection{}

Let $Y$  be a smooth projective curve of genus $g\geq 2$ over a number field
$K$. The $L$-function of $Y$ is an analytic function of one complex variable
$s$ defined for $\Re(s)\geq 2$ as an Euler product
\[
   L(Y,s) := \prod_\p L_\p(Y,s),
\]
where $\p$ ranges over 
the prime ideals of $K$  and the local $L$-factor is of
the form
\begin{equation}\label{eqn:poly_p}
        L_\p(Y,s) = \frac1{P(({\rm N}\p)^{-s})}.
\end{equation}
Here ${\rm N}\p$ is the norm of $\p$ and $P(T)=1+\ldots\in\ZZ[T]$ is a
polynomial with integer coefficients depending on $\p$. Another invariant
associated with $Y$ and relevant for our discussion is the {\em conductor of
  the $L$-function}. It is a positive real number of the form
\[
        N:= \delta_K^{2g}\cdot\prod_\p({\rm N}\p)^{f_\p},
\]
where $\delta_K$ is the discriminant of $K$ and $f_\p\geq 0$ is a nonnegative
integer called the {\em conductor exponent} at $\p$, which is zero for almost
all $\p$. We refer to \cite{superell}, \S~2, for precise definitions of
$L_\p(Y,s)$ and $f_\p$. 

It is conjectured that $L(Y,s)$ has an analytic continuation to the whole
complex plane, and a functional equation of the form
\begin{equation} \label{eqn:funceq}
      \Lambda(Y,s) = \pm\Lambda(Y,2-s),
\end{equation} 
where 
\[
   \Lambda(Y,s):=N^{s/2}(2\pi)^{-gs}\Gamma(s)^gL(Y,s).
\]

\subsection{}

The main motivation for this paper -- which continues a project begun in
\cite{superell} -- is the question how to compute the local $L$-factor
$L_\p(Y,s)$ and the conductor exponent $f_\p$ explicitly, given the curve $Y$
and the prime $\p$. If the curve $Y$ has good reduction at $\p$ (which is true
for almost all $\p$) it is well known how to do that. Namely, $f_\p=0$ and
$L_\p(Y,s)=P(\Yb,({\rm N}\p)^{-s})^{-1}$, where $P(\Yb,T)\in\ZZ[T]$ is the
numerator of the zeta function of the reduction $\Yb$ of $Y$ at $\p$. To
compute $P(\Yb,T)$ for small primes $\p$, one can simply count the number of
$\FF_{q^n}$-rational points on $\Yb$, for $q={\rm N}\p$ and
$n=1,\ldots,g$. The complexity of this approach is bounded by
$\OO(q^g)$. There is an extensive literature dealing with various methods for
lowering this asymptotic bound, see e.g.\ \cite{GaudryHarley2000},
\cite{Kedlaya2001} or \cite{KedlayaSutherland}.

If $Y$ has bad reduction at $\p$ then it is not so easy to compute $L_\p(Y,s)$
and $f_\p$ directly from the curve $Y$, even if ${\rm N}\p$ is very small
(say, equal to $p=2$). Prior to \cite{superell}, no general and systematic
approach was known except for $g=1$, and for $g=2$ and $q$ odd. Nevertheless,
there have been successful attempts to compute $L_\p(Y,s)$ and $f_\p$ for all
$\p$ for curves of genus $g\geq 2$, for instance by Dokchitser, de Jeu and
Zagier (\cite{DokdeJeuZagier}), and by Booker (\cite{booker},
\cite{booker05}). Without going into details, their method is either
based on guessing $L_\p(Y,s)$ and $f_\p$ for the finitely many primes of bad
reduction and then verifying this guess by checking the functional equation
\eqref{eqn:funceq}, or by tailoring $L_\p(Y,s)$ and $f_\p$ in such a way that
\eqref{eqn:funceq} holds. Of course, as long as \eqref{eqn:funceq} remains a
conjecture, these methods are unable to prove correctness of the result.

In contrast, in this paper we compute $L_\p(Y,s)$ and $f_p$ directly
for all primes $\p$ of bounded size (i.e.\ for ${\rm N}\p\leq M$ for a
certain constant $M$), for many curves $Y$ over $\QQ$ of genus
$g=2,3,4,5,6$. The computation at the primes of bad reduction is done
using the methods of \cite{superell} and \cite{superp}, and they are
provably correct. We then verify the functional equation
\eqref{eqn:funceq} numerically, using the methods of
\cite{Dokchitser04}. 

The expected sign of the functional equation \eqref{eqn:funceq} is
also known to be a product of local factors, the so called {\em local
  root numbers} (\cite{Deligne}).  In principal, it should be possible
to compute the local root numbers using our methods, but we have not
tried to do that. Of course, as a side effect of our numerical
verification of the functional equation we obtain an experimental
value for its sign, which is correct with
very high probability.

We expect that with our approach it is now possible to compute examples with
much larger conductor $N$ and genus $g$ than before. However, we have not
tried to push computations to their limit. In our largest example, $g=6$ and
$N=7 \cdot 11 \cdot 13 \cdot 89 \cdot 431 \cdot 857 \approx 3 \cdot 10^{10}$,
which is comparable to the largest examples considered in
\cite{DokdeJeuZagier}.

\subsection{}

In \cite{superell} we have shown how to compute the local factor $L_\p(Y,s)$
and the conductor exponent $f_\p$ from the {\em stable reduction} of $Y$ at
$\p$. Furthermore, we have shown how this can be done explicitly for {\em
  superelliptic curves}, i.e.\ curves $Y$ given by an equation of the form
\[
    y^n= f(x),
\]
where $f\in K[x]$ is a polynomial with coefficients in $K$. A
serious restriction that we  imposed in \cite{superell} is that the
exponent $n$ is prime to the residue characteristic of $\p$. This
restriction can be removed, using the results of \cite{ArzdorfDiss} and
\cite{superp}. So in principal we can compute $L_\p(Y,s)$ and $f_\p$ for all
primes $\p$ and all superelliptic curves. There is also no fundamental
difficulty to extend our methods to curves which are not
superelliptic. However, the details can get tricky, and it is rather hard
to implement algorithms which work for general classes of curves. 

The main class of examples we consider in the present paper is constructed in
a way to illustrate our main point, while being at the same time as simple as
possible and to be manageable by a straightforward algorithm. We consider a
rather general family of hyperelliptic curves over $\QQ$ of fixed genus
$g_Y\geq 2$. Within this family we search for examples of curves $Y$ which
have semistable reduction at every prime number $p$. For each curve
satisfying this condition we compute its $L$-series and conductor and
numerically verify the functional equation. 

If a curve $Y$ does not have semistable reduction at a prime $\p$, but
only after replacing $K$ by a finite extension, the computation of
$L_\p(Y,s)$ and $f_\p$ is much more involved. At the moment, we have not
yet implemented algorithms which can handle such examples in a
routine fashion. We discuss three examples in  detail,
illustrating the difficulties occurring.  The discussion illustrates
that each individual problem can typically be solved by a
knowledgeable human supported by customized computational tools.

\subsection{}

The structure of the paper is as follows. In \S~\ref{sec:cohom} we recall how
to compute the local $L$-factor and the conductor exponent at a prime $\p$
where the curve $Y$ has semistable reduction. The explicit expression for the
local $L$-factor (resp.\ the conductor exponent) can be found in Proposition
\ref{prop:badfactor} (resp.\ Corollary \ref{cor:badfactor}).

In \S~\ref{sec:ssred} we consider a rather general class of
hyperelliptic curves, and determine necessary and sufficient conditions
for these curves to have semistable reduction everywhere (Lemmas
\ref{lem:redp=2} and \ref{lem:podd}).  \S~\ref{sec:algo} summarizes
the algorithm for computing the local $L$-factor and the conductor
exponent at the primes of bad reduction of the curves satisfying these
conditions, and for verifying the functional equation numerically.
Examples are given in \S~\ref{sec:exa}. In \S~\ref{sec:more_examples}
three examples of superelliptic curves which do not have semistable
reduction everywhere are discussed.

All data from the examples discussed in this paper can be retrieved from\\ 
{\tt https://www.uni-ulm.de/index.php?id=64504}

\section{\'Etale cohomology of a semistable curve}\label{sec:cohom}

Let $Y$ be a smooth projective and absolutely irreducible curve of genus
greater than or equal to $2$ defined over a number field $K$. In this section
we recall from \S~2 of \cite{superell} the description of the local $L$-factor
and the conductor exponent at a prime $\p$ of $K$ in the case that $Y$ has
semistable reduction at $\p$. In general, the curve $Y$ only admits semistable
reduction after passing to a finite extension. The main result of this section
gives an explicit, computable expression for the local $L$-factor and the
conductor exponent in the case that no field extension is needed.

\subsection{}

For a (finite) prime $\p$ of $K$ we write ${\mathcal O}_\p\subset K$ for the
local ring and $\FF=\FF_\p$ for the residue field. Let $q:=\rm N\p=\abs{\FF}$
denote the norm of $\p$.

Throughout this section we assume that $Y$ has \textsl{semistable reduction}
at $\p$.  Recall that this means that there exists a proper and flat model
$\Y$ of $Y$ over $\mathcal O_\p$ whose special fiber
$\bar{Y}:=\Y\otimes_{\OO_\p} \FF$ is semistable, i.e.\ $\bar{Y}$ is reduced
and has only ordinary double points as singularities. We keep the semistable
model $\Y$ fixed and call its special fiber $\bar{Y}$ the \textsl{semistable
  reduction} of $Y$ at $\p$ (even though $\Yb$ is not uniquely determined
without further assumptions). We write $\bar{Y}_k:=\bar{Y}\otimes_{\FF} k$ for
the base change of $\bar{Y}$ to the algebraic closure $k$ of $\FF$. We
denote the absolute Galois group of $\FF$ by $\Gamma_{\FF}$. Let
$\Frob_\p\in\Gamma_\FF$ denote the arithmetic Frobenius element, i.e.\ the
element determined by
\[
       \Frob_\p(a)=a^q
\]
for $a\in k$.

If $\Yb$ is smooth (i.e.\ $\p$ is a prime of good reduction) then it is well
known that the local $L$-factor $L_\p(Y,s)$ may be computed by point
counting on $\bar{Y}$. Moreover, the conductor exponent is zero. In our case
(where $\Yb$ is semistable) this generalizes as follows. Let
$H^i_\et(\bar{Y}_k, \QQ_\ell)$ be the $i$th $\ell$-adic \'etale cohomology
group of $\bar{Y}_k$, where $\ell$ is an auxiliary
prime different from the residue characteristic of $\p$. Write
$\Frob_{\Yb}:\bar{Y}\to \bar{Y}$ for the relative $\FF_\p$-Frobenius
morphism. For $n\in\NN$ let $\FF_n\subset k$ be the (unique) finite extension
of $\FF$ of degree $n$. The {\em zeta function} of $\Yb$ is defined as
\[
   Z(\Yb,T):=\exp\Big(\sum_{n\geq 1} \abs{\Yb(\FF_n)}\cdot\frac{T^n}{n}\Big).
\] 
It is well known that $Z(\Yb,T)$ is a rational function of the form
\[
    Z(\Yb,T) = \frac{P_1(T)}{P_0(T)\cdot P_2(T)},
\]
where 
\[
     P_i(T):=\det(1-T\cdot\Frob_{\Yb}\mid_{H^i_\et(\bar{Y}_k, \QQ_\ell)}).
\]
See e.g.\ \cite{MilneEC}, Theorem 13.1. 

\begin{prop} \label{prop:Lfactor} 
  The local $L$-factor is given by the formula
  \[
     L_\p(Y/K, s)=P_1(q^{-s})^{-1},
  \]
  where $P_1(T)\in\ZZ[T]$ is the numerator of the zeta function of $\Yb$. The
  conductor exponent is 
  \[
      f_\p=2g_Y-\dim H^1_\et(\bar{Y}_k, \QQ_\ell) = 2g_Y-\deg(P_1).
  \]
\end{prop}

\proof This follows directly from \cite{superell}, Corollaries 2.5 and
2.6, since we assume that $Y$ has semistable reduction over $K$.  \Endproof

\begin{rem} \label{rem:pointcount} Assuming we have an explicit equation for
  the curve $\Yb$, Proposition \ref{prop:Lfactor} gives a simple way of
  computing $L_\p(Y/K,s)$ via point counting. Indeed, it suffices
  to compute the polynomials $P_i$ for $i=0,1,2$. For $i=0,2$ this is
  easy. Since $\Yb_k$ is connected we have $H^0_\et(\Yb_k,\QQ_\ell)=\QQ_\ell$,
  with trivial action of $\Frob_{\Yb}$ and hence
  \[
     P_0(T) = 1-T.
  \]
  Let $\Yb_i$ denote the irreducible components of $\Yb$, and let $m_i$
    denote
  the number of irreducible components of $\Yb_i\otimes k$. Then
  \[
     P_2(T) = \prod_i \big(1-(qT)^{m_i}\big).
  \] 
  So in order to compute $P_1(T)$ (which has degree $\leq 2g_Y$) it suffices
  to know the first $2g_Y+1$ terms of the power series $Z(\Yb,T)$, which 
  depend on $\abs{\Yb(\FF_n)}$ for $n=1,\ldots,2g_Y$. 
\end{rem}

\begin{rem} \label{rem:pointcount2}
  If $\Yb$ is smooth, the bound $2g_Y$ from Remark
  \ref{rem:pointcount} can be improved to $g_Y$, using the functional
  equation.  More precisely, the polynomial $P_1(T)$ has the form
  \[
         P_1(T) = c_0+c_1T+\ldots+c_{2g_Y}T^{2g_Y} \in\ZZ[T],
\]
with $c_0=1$ and satisfies the functional equation
% \[
%      P_1(1/qT) = q^{1-g_Y}T^{2-2g_Y}P_1(T),
% \]
% \michel{das ist die F-Eq fuer die Zetafunktion!}
% \michel{richtig ist:}
\[
     P_1(1/qT) = q^{-g_Y}T^{-2g_Y}P_1(T),
\]
see e.g.\ \cite{MilneEC}, Theorem 12.6. The functional equation is equivalent
to%\michel{exponent of $q$ corrected}
\[
     c_{2g_Y-i} = q^{g_Y-i}c_i, \quad i=0,\ldots,g_Y.
\]
This means that $P_1$ is already determined by the coefficients
$c_0,\ldots,c_{g_Y}$. It follows that $P_1(T)$ can
be computed by counting the number of points of $\Yb$ over the fields
$\FF_n$, for $n=1,\ldots,g_Y$. 
%\irene{Ich habe $g_0$ hier wieder durch $g_Y$ ersetzt.}
\end{rem}

\subsection{}

\label{sec:cohomsub2}
We have seen in the previous section that we can compute the local $L$-factor
and the conductor exponent of $Y$ at a prime $\p$ of semistable reduction,
provided we know an explicit equation for the stable reduction $\Yb$. In order
to do this, we have to count the number of points of $\Yb$ over certain finite
extensions of the residue field of $\p$. We note in passing that all the
computations done for the present paper only use the naive counting method (as
opposed to more sophisticated methods as e.g.\ in \cite{Kedlaya2001} or
\cite{GaudryHarley2000}).

If $\Yb$ is smooth (i.e.\ if $Y$ has good reduction at $\p)$ we
can use the functional equation to reduce the cost of point counting
drastically (Remark \ref{rem:pointcount2}).  In this section we
 extend this trick to the case where $\Yb$ is semistable. To
keep the discussion simple, we assume that the curve $\Yb$ is
absolutely irreducible. This assumption is satisfied for our main
class of examples considered in \S~\ref{sec:ssred}. As a first
consequence we see that the denominator of the zeta function is of the
most simple form,
\[
     Z(\Yb,T) = \frac{P(\Yb,T)}{(1-T)(1-qT)}.
\]
Here $P(\Yb,T)=P_1(T)$ in the notation of the previous subsection. 

Let
\[
    \pi:\Yb_0\to\Yb
\]
denote the normalization of $\Yb$. Then $\Yb_0$ is a smooth projective
absolutely irreducible curve and $\pi$ is a finite  birational morphism. If
$\xi\in\Yb$ is a closed point then the fiber $\pi^{-1}(\xi)$ has degree one
over $\FF(\xi)$ if $\xi$ is a smooth point and  has degree two if $\xi$ is an
ordinary double point. In the latter case, we say that $\xi$ is a {\em split}
(resp.\ a {\em nonsplit}) double point if $\pi^{-1}(\xi)$ consists of two points
(resp.\ of one point). Geometrically the map $\pi$ may be visualized as in
Figure \ref{fig:ybar_y0bar}. 

\begin{figure}[h!]\centering 
\begin{tikzpicture}[scale = 0.12]
 \coordinate (A) at (-4,0);
 \coordinate (B) at  (0,0);
 \coordinate (C) at  (15,11);
 \coordinate (C1) at  (5,11);
 \coordinate (D) at  (20,1);
 \coordinate (E) at  (0,20);
 \coordinate (F) at  (35,12);
 \coordinate (F1) at  (25,12);
 \coordinate (G) at  (40,4);
 \coordinate (H) at  (11,25);
 \coordinate (H1) at  (38,22);
 %\coordinate (H2) at  (50,20);
 \coordinate (H3) at  (60,28);
 \coordinate (K) at  (55,13);
 \coordinate (L) at  (45,13);
 \coordinate (M) at  (60,3);

 \coordinate (I) at (0,-10);
 \coordinate (J) at (60,-10);

  \coordinate (s1) at (9.02,5.68);
  \coordinate (s10) at (9.02,-10);
   \coordinate (s2) at (31.0,6.39);
  \coordinate (s20) at (31.0,-10);
   \coordinate (s3) at (50.85,8.25);
  \coordinate (s30) at (50.85,-10);

%\draw [gray!55]   (B) -- (C) -- (C1) -- (D)  -- (F) -- (F1) -- (G) --(K)--(L)--(M) ;

\draw [black] plot [smooth, tension=1.0] coordinates {   (B) (C) (C1) (D)   (F) (F1) (G)  (K) (L)(M)} ;

\draw [black] plot [smooth, tension=1.0] coordinates {  (E)(H)(H1)(H3)} ;

\node [right] at (62,28) {$\bar Y_0$};

\node [right] at (62,3) {$\bar Y$};
\draw [black ]     (I) -- (J);
\node [right] at (62,-10) {$\bar X$};

\draw [black , dashed]     (s1) - - (s10);

 \filldraw [black]      (s1) circle (5pt) ;
 \filldraw [black]      (s10) circle (5pt) ;

\draw [black , dashed]     (s2) - - (s20);

 \filldraw [black]      (s2) circle (5pt) ;
 \filldraw [black]      (s20) circle (5pt) ;

\draw [black , dashed]     (s3) - - (s30);

 \filldraw [black]      (s3) circle (5pt) ;
 \filldraw [black]      (s30) circle (5pt) ;

\foreach \i in {-3,...,3}
{
\draw [white, line width = 1.2 pt]  (40+\i/1.5,-11) -- (40+\i/1.5,30);

}

\end{tikzpicture}
{\caption{Normalization of $\Yb$}\label{fig:ybar_y0bar}}
\end{figure}

Let $g_0$ denote the genus of $\Yb_0$. The zeta function of $\Yb_0$ has the
form
\[
     Z(\Yb_0,T) = \frac{P(\Yb_0,T)}{(1-T)(1-qT)},
\]
where $P(\Yb_0,T)$ satisfies the functional equation,
and hence can be determined by counting $\abs{\Yb_0(\FF_n)}$ for
$n=1,\ldots,g_0$ (Remark \ref{rem:pointcount2}).

The following result reduces the calculation of the local $L$-factor
in our situation to point counting on the normalization $\bar{Y}_0$ of
$\bar{Y}$. 

\begin{prop}\label{prop:badfactor} 
  Let $\Sc$ denote the set of singular points of $\Yb$. For $\xi\in\Sc$ we let
  $d_\xi:=[\FF(\xi):\FF]$ denote the degree of $\xi$. Furthermore, define
  $\varepsilon_\xi:=1$ (resp.\ $\varepsilon_\xi:=-1$) if $\xi$ is a split (resp.\ a
  nonsplit) double point. Then 
  \[
     P(\bar{Y}, T)=P(\bar{Y}_0, T)\cdot 
           \prod_{\xi\in \Sc}(1-\varepsilon_\xi T^{d_\xi}).
  \]
\end{prop}

\proof Lemma 2.7.(1) of \cite{superell} implies that the $\ell$-adic
\'etale cohomology group of $\bar{Y}$ decomposes as a direct sum of
$\Gamma_\FF$-modules
\[
    H^1_\et(\bar{Y}_k,\QQ_\ell)=H^1_\et(\bar{Y}_{0,k},\QQ_\ell)
               \oplus H^1(\Delta_{\bar{Y}_k},\QQ_\ell),
\]
where $\Delta_{\bar{Y}_k}$ denotes the  graph of components of $\Yb_k$. Therefore,
it suffices to show that 
\[
    \det(1-T\cdot\Frob_\FF|_{H^1(\Delta_{\bar{Y}_k}, \QQ_\ell)}) = 
      \prod_{\xi\in \Sc}(1-\varepsilon_\xi T^{d_\xi}).
\]
This amounts to computing the character of the representation of
$\Gamma_\FF$ acting on $H^1(\Delta_{\bar{Y}_k}, \QQ_\ell)$, which is
described in Lemma 2.7.(2) of \cite{superell}.

Since we assume that $\Yb$ is a semistable, absolutely irreducible
curve, the graph $\Delta_{\Yb_{k}}$ is a bouquet of $r$ circles, where
\[
     r = \sum_{\xi\in \Sc} d_\xi
\]
is the number of ordinary double points of $\Yb_{k}$ (see Figure
\ref{fig:ybar_y0bar}). An element $\xi\in \Sc$ corresponds to a
$\Gamma_\FF$-orbit of edges of $\Delta_{\Yb_{k}}$. Furthermore, $\xi$ is a
split (resp.\ nonsplit) ordinary double point if and only if the stabilizer
$\Gamma_{\FF(\xi)}$ acts trivially (resp.\ acts by reversing orientation) on
any of the edges in the orbit corresponding to $\xi$. Lemma 2.7.(2) of
\cite{superell} implies that the character of $H^1(\Delta_{\bar{Y}_k},
\QQ_\ell)$ considered as $\Gamma_\FF$-representation is
\begin{equation}\label{eqn:char}
  \chi_{\rm sing}:=\bigoplus\limits_{\xi\in \Sc} 
          \Ind_{\Gamma_{\FF(\xi)}}^{\Gamma_{\FF}}\varepsilon_\xi.
\end{equation}
Here we interpret the integer $\varepsilon_\xi\in \{\pm1\}$ as the character of
a $1$-dimensional representation of $\Gamma_{\FF(\xi)}$. Namely,
$\varepsilon_\xi$ is the trivial character if $\varepsilon_\xi=1$ and the
unique character of order $2$ if $\varepsilon_\xi=-1$. 
The statement of the proposition now follows from an elementary
calculation. 

For a proof which does not use \'etale cohomology, see \cite{AubryPerret}.
\Endproof

\begin{cor} \label{cor:badfactor}
  The conductor exponent is
  \[
     f_\p = r = \sum_{\xi\in \Sc} d_\xi.
  \]
\end{cor}

\section{Hyperelliptic curves with semistable reduction 
           everywhere}
\label{sec:ssred}

In this section we consider a class of hyperelliptic curves of genus greater
than or equal to $2$ which are defined over a number field $K$.  We find
conditions on the equation which guarantee that the curve has semistable
reduction at every prime. This makes it relatively easy to calculate the local
$L$-factor at the bad primes, even for residue characteristic $p=2$.

\subsection{} \label{ssred1}

We fix a number field $K$, an integer $g_Y\geq 2$ and two polynomials
$g,h\in\OO_K[x]$ satisfying the following three conditions:
\begin{itemize}
\item The polynomial $g$ is monic of degree $2g_Y+1$.
\item
  The degree of $h$ is at most $g_Y$. 
\item
  The polynomial $f:=4g+h^2$ has no multiple roots. 
\end{itemize}
Let $Y$ be the smooth projective curve over $K$ which is given birationally by
the equation
\begin{equation}\label{eqn:p=2}
       y^2 + h(x)y = g(x).
\end{equation}
By this we mean that the function field of $Y$ is the field extension of $K$
with two generators $x,y$ satisfying the above equation. 
Our assumptions imply that $Y$ is absolutely irreducible and, more
specifically, a hyperelliptic curve of genus $g_Y$. An alternative equation
for $Y$ is
\begin{equation} \label{eqn:oddprime}
     u^2 = f(x):=4g(x)+h(x)^2, 
\end{equation}
where $u:=2y+h(x)$. Depending on the residue characteristic considered, either
(\ref{eqn:p=2}) or (\ref{eqn:oddprime}) will be more useful.

Equation (\ref{eqn:p=2}) defines a smooth plane curve with a unique point `at
infinity' which we denote by $\infty$.  It will be useful to
have a similar equation for a neighborhood of the point $\infty$. For this we
set $z:=x^{-1}$, $w:=z^{g_Y+1}y$, $g^*:=z^{2g_Y+2}g$ and
$h^*:=z^{g_Y+1}h$. Considering $g^*,h^*$ as polynomials in $z$, \eqref{eqn:p=2}
can be rewritten as 
\begin{equation} \label{eqn:infty}
    w^2+h^*(z)w = g^*(z).
\end{equation}
This is again an equation for a smooth plane curve, and the point $\infty$ has
coordinates $(z,w)=(0,0)$. Note that we have used the assumption that $g$ has
odd degree to prove smoothness at $\infty$. 

\subsection{} \label{ssred2}

We now choose a prime ideal $\p$ of $\OO_K$. Let $\OO_\p$ denote the local
ring and $\FF_\p$ the residue field of $\p$, as in \S~\ref{sec:cohom}.

Let $\X=\PP^1_{\OO_\p, x}$ be the projective line over $\OO_\p$ with
coordinate $x$ and write $\Y$ for the normalization of $\X$ in the function
field $K(Y)$ of $Y$. Then $\Y$ is a {\em model} of $Y$ over $\OO_\p$, i.e.\
$\Y$ is a flat and proper $\OO_\p$-scheme of finite type with generic fiber
$Y$.

We denote by $\Yb$ and $\Xb$ the special fibers of $\Y$ and $\X$,
respectively. These are proper curves over $\FF_\p$, and
$\Xb=\PP^1_{\FF_\p}$. We write $\xb,\yb$ for the image of $x,y$ in the
function ring of $\Yb$, and $\bar{g}$ (resp.\ $\bar{h}$) for the image of $g$
(resp.\ $h$) in $\FF_\p[\xb]$. The following proposition shows that the curve
$\bar{Y}$ is completely determined by the `reduction' of (\ref{eqn:p=2})
modulo $\p$.

\begin{prop}\label{prop:givenby} 
  The curve $\Yb$ is reduced and absolutely irreducible. The point $\infty$
  reduces to a smooth point $\bar{\infty}\in\Yb$, and the affine open part
  $\Yb-\{\bar{\infty}\}$ is a plane curve with equation
  \begin{equation} \label{eqn:Ybeq}
     \yb^2+\hb(\xb)\yb = \gb(\xb).
  \end{equation}
\end{prop}

\proof 
The curve $\X$ has an open affine covering $\{\Spec A_1,\Spec A_2\}$, where
$A_1=\OO_\p[x]$ and $A_2=\OO_\p[z]$. It follows that $\Y$ has an open affine
covering $\{\Spec B_1,\Spec B_2\}$, where $B_i$ is the integral closure of
$A_i$ in $K(Y)=K(x,y)$. We claim that $B_1=A_1[y]$ and $B_2=A_2[w]$. 

Let us first consider $B_1$. The minimal polynomial for $y$ over the function
field $K(X)=K(x)$ is a monic polynomial with coefficients in $A_1$, 
\[
        F_1:=T^2+hT-g.
\]
It follows that $B_1':=A_1[y]\cong A_1[T]/(F_1)$ is finite and flat over
$A_1$. Moreover, $B_1'\otimes_{\OO_\p}K$ is integrally closed because
\eqref{eqn:p=2} defines a smooth curve. Now
\[
     B_1\otimes\FF_\p = \FF_\p[x,T]/(\bar{F}_1),
\]
where 
\[
     \bar{F}_1=T^2+\hb T-\gb
\]
is the image of $F_1$ in $\FF_\p[x,T]$. The polynomial $\gb\in\FF_\p[x]$ still
has odd degree $2g_Y+1$ (because we have assumed that $g$ is monic). It
follows that $\bar{F}_1$ is absolutely irreducible. We conclude using
Lemma 4.1.18 of \cite{liu} that $B_1'$ is integrally closed and hence
$B_1=A_1[y]$.  The proof that $B_2=A_2[w]$ is similar; one uses that
\eqref{eqn:infty} defines a smooth plane curve which remains reduced and
irreducible after reduction to the residue field. The remaining statements are
also easy to show.
\Endproof

\subsection{}\label{sec:redp=2}

We continue with the notation and assumptions of \S\S~\ref{ssred1} and
\ref{ssred2}. Additionally, we assume that the residue field $\FF_\p$
has characteristic $p=2$.

\begin{notation} \label{not:xi} Let $\xi\in\Yb-\{\bar{\infty}\}$ be a closed
  point and $\FF_p(\xi)$ the residue field of $\xi$. We consider
  $\Yb-\{\bar{\infty}\}$ as an affine plane curve with coordinate functions
  $\xb,\yb$. Set $a:=\xb(\xi),b:=\yb(\xi)\in\FF_p(\xi)$. Then
  $\FF_\p(\xi)=\FF_\p(a,b)$, and we write $\xi=(a,b)$. 
\end{notation}

\begin{lem} \label{lem:redp=2}
  Let $\xi=(a,b)\in\Yb-\{\bar{\infty}\}$  be a closed point. 
  \begin{enumerate}
  \item
    The point $\xi$ is a singularity of $\Yb$ if and only if
    \begin{equation} \label{eqn:redp=2.0}
         \hb(a)=0=(\hb'(a))^2\gb(a)+(\gb'(a))^2.
    \end{equation}
    Here $\hb',\gb'\in\FF_\p[\xb]$ are the formal derivatives of $\hb,\gb$ with
    respect to $\xb$.  
  \item
    Assume that $\xi$ is a singularity. Then $\xi$ is an ordinary double point
    if and only if $\hb'(a)\neq 0$.
  \end{enumerate}
\end{lem}

\proof
Clearly, $\xi=(a,b)$ satisfies \eqref{eqn:Ybeq}:
\begin{equation} \label{eqn:redp=2.1}
  b^2+\hb(a)b=\gb(a).
\end{equation}
The Jacobian criterion says that $\xi$ is singular if and only if 
\begin{equation} \label{eqn:redp=2.2}
     \hb(a)=0,\quad \hb'(a)b=\gb'(a).
\end{equation}
Using \eqref{eqn:redp=2.1} to eliminate $b$, we see that
\eqref{eqn:redp=2.1} is equivalent to \eqref{eqn:redp=2.0}. Now (i) is
proved.

For the proof of (ii) we assume that $\xi$ is singular and compute the tangent
cone of $\Yb$ at $\xi$, using \eqref{eqn:redp=2.1} and
\eqref{eqn:redp=2.2}. We obtain
\[
    (\yb+b)^2 + \hb'(a)(\yb+b)(\xb+a) + \gb_2(\xb+a)^2 = 0,
\]
where $\gb_2$ is the coefficient of $\xb^2$ in the Taylor expansion of $\gb$
at $\xb=a$. As we are in characteristic $2$, the underlying quadratic form is
nondegenerate if and only if $\hb'(a)\neq 0$. This proves (ii).
\Endproof

\begin{cor}\label{cor:redp=2} 
  The curve $\Yb$ is semistable if and only if
  $\hb\neq 0$ and
  \[
      \gcd(\bar h, \bar h',\bar g')=1.
  \]
\end{cor}

\proof To prove the corollary it suffices to show that
$\xi=(a,b)\in\Yb-\{\bar{\infty}\}$ is a smooth or   an ordinary double point
if and only if $(\hb(a),\hb'(a),\gb'(a))\neq (0,0,0)$. This follows directly
from Lemma \ref{lem:redp=2}.
\Endproof

From now on we assume that $\Yb$ is semistable, and we  use the results
from \S ~\ref{sec:cohom} to compute the local $L$-factor and the conductor
exponent of $Y$ at $\p$. Let
\[
       \pi:\Yb_0\to\Yb
\]
be the normalization of $\Yb$. Recall that $\pi$ is a finite map which is an
isomorphism above the smooth locus of $\Yb$ (as in Figure
\ref{fig:ybar_y0bar}). In order to use Proposition \ref{prop:badfactor} and
Corollary \ref{cor:badfactor} we need to know the set $\Sc$ of singular points
of $\Yb$, the invariants $d_\xi$ and $\varepsilon_\xi$, for all $\xi\in\Sc$, and
an explicit equation for $\Yb_0$. This will be achieved by the following
proposition and Corollary \ref{cor:redp=2.2}.

\begin{prop} \label{prop:redp=2}
  Assume that $\Yb$ is semistable. Set
  \[
       r:=\gcd(\hb,(\hb')^2\gb+(\gb')^2)\in\FF[\xb].
       \]
             Then the following holds.
  \begin{enumerate}
  \item
    A point $\xi\in\Yb-\{\bar{\infty}\}$ is singular if and only if $r(a)=0$. 
  \item
    The polynomial $r$ is separable, i.e.~all roots of $r$ over the
    algebraic closure $k$ of $\FF$ are simple. Furthermore,
    $\tilde{h}:=\hb/r\in\FF[\xb]$ is prime to $r$.
  \item
    There exists $s\in\FF[\xb]$   such that
    \[
         r^2 \mid \gb+s^2+\hb s.
    \]
  \item
    Set $\tilde{g}:=(\gb+s^2+\hb s)/r^2\in\FF[\xb]$ and
    $\tilde{y}:=(\yb+s)/r\in\FF(\Yb)$. Then $\tilde{y}$ satisfies 
    \begin{equation} \label{eqn:redp=2.3}
      \tilde{y}^2+\tilde{h}\tilde{y} = \tilde{g},
    \end{equation}
    which is an equation for the smooth plane affine
    $\Yb_0-\{\bar{\infty}\}$.
  \end{enumerate}
\end{prop}

\proof (i) follows directly from Lemma \ref{lem:redp=2} (i). Now assume that
$a$ is a root of $r$. Then there is a unique point $\xi=(a,b)\in\Yb$, and it
is a singularity. Since we assume that $\Yb$ is semistable, $\xi$ is even an
ordinary double point. Therefore, it follows from Lemma \ref{lem:redp=2} (ii)
that $\hb'(a)\neq 0$. We conclude that all roots of $r$ are simple roots of
$\hb$. This proves (ii).

Since $\FF$ is a perfect field of characteristic $2$ and $r$ is
separable by (ii), there exists a polynomial $s\in\FF[\xb]$ such that
\[
       s^2\equiv \gb \pmod{r}.
\]
Set $\tilde{y}:=(\yb+s)/r\in\FF(\Yb)$. Then $\tilde{y}$ satisfies equation
\eqref{eqn:redp=2.3}. 

For the proof of (iii) we have to show that $\tilde{g}$ is a
polynomial. Assume that $a\in k$ is a pole of $\tilde{g}$. By (ii)
 $r$ has a simple zero at $a$. The
choice of $s$ implies that $\hb$ also has a simple zero at $a$, and
hence that $\tilde{g}$ has a simple pole at $\xb=a$. But this would
mean that the map $\Yb_0\to\Xb=\PP^1_\FF$ is branched at $\xb=a$.
This would imply that there exists a unique smooth point
$\xi=(a,b)\in\Yb$ above $\xb=a$, contradicting the fact that
$r(a)=0$. Now (iii) is proved.
         
It follows from (iii) that there is a finite birational morphism $\Yb_1\to\Yb$
which is an isomorphism at $\bar{\infty}$ and such that
$\Yb_1-\{\bar{\infty}\}$ is the plane affine curve given by
\eqref{eqn:redp=2.3}. Let $\xi=(a,b)\in\Yb_1-\{\bar{\infty}\}$ be a closed
point. If $\xi$ is a singular point, then $\tilde{h}(a)=0$ by the Jacobian
criterion. But then $r(a)\neq 0$ by (ii) and the definition of
$\tilde{h}$. Therefore, $\xi$ lies above a smooth point of $\Yb$. Since
$\Yb_1\to\Yb$ is finite, it follows that $\xi$ is a smooth point as well,
contradiction. We conclude that $\Yb_1$ is smooth. This implies that
$\Yb_1=\Yb_0$ is the normalization of $\Yb$ and completes the proof of
the proposition.
\Endproof 

\begin{cor} \label{cor:redp=2.2}
   Assume that $\Yb$ is semistable. 
  \begin{enumerate}
  \item
    There is a bijection between the set $\Sc$ of singular points of $\Yb$ and
    the irreducible factors of the polynomial $r\in\FF[\xb]$ defined in
    Proposition \ref{prop:redp=2}.  
  \item
    A singular point $\xi=(a,b)\in\Yb$ is a split (resp.\ a non split)
    ordinary double point if the polynomial
    \[
         T^2+\tilde{h}(a)T+\tilde{g}(a) \in\FF[T]
    \]
    is reducible (resp.\ irreducible).
  \item
    The conductor exponent at $\p$ is 
    \[
        f_\p = \deg(r).
    \]
  \end{enumerate}
\end{cor}

\subsection{}\label{sec:redpodd}

We now switch to the case of a prime $\p$ with residue characteristic $p\geq
3$. It will be more convenient to use Equation $(\ref{eqn:oddprime})$ to
describe the curve $Y$: 
\[
    u^2 = f(x):=4g(x)+h(x)^2.
\] 
Recall that this equation is derived from \eqref{eqn:p=2} by the substitution
$y=(u-h)/2$. Since $2$ is a unit in $\OO_\p$, the same substitution works for
the model $\Y$. It follows that the special fiber $\Yb$ of $\Y$ is given by
the equation
\begin{equation} \label{eqn:podd1}
    \ub^2 = \fb(\xb).
\end{equation}
Here $\fb\in\FF_\p[\xb]$ denotes the image of $f$ in $\FF_\p[\xb]$ and $\ub$
the image of $u$ in $\FF(\Yb)$. We also adopt Notation \ref{not:xi} to this
new equation and write a closed point $\xi\in\Yb-\{\bar{\infty}\}$ in the form
$\xi=(a,b)$, where $(a,b)$ is a solution to \eqref{eqn:podd1}. 

Note that by choice of $g$ (monic, degree $2g_Y+1$) and $h$ (degree $\leq
g_Y$), both $f$ and $\bar f$ have degree $2g_Y+1$. The polynomial $f$ is
separable by assumption, but in general this will not be true for $\fb$. The
following is very easy to show.

\begin{lem} \label{lem:podd}
  The curve $\Yb$ is semistable if and only if $\fb$ has at most double
  roots. 
\end{lem}

Let us assume from now on that the curve $\Yb$ is semistable. By Lemma
\ref{lem:podd}, the polynomial $\fb$ has at most double roots. It follows
that there is a unique decomposition
\[
     \fb = r^2\cdot s,
\]
where $r,s\in \FF[\xb]$ are separable and relatively prime. 

\begin{prop} \label{prop:podd}
  We assume that $\Yb$ is semistable. Let $\xi=(a,b)\in\Yb-\{\bar{\infty}\}$
  be a closed point. 
  \begin{enumerate}
  \item
    The point $\xi$ is a singularity of $\Yb$ if and only if $r(a)=0$.
  \item
    Assume $\xi$ is a singularity. Then $\xi$ is a split (resp.\ a non split)
    ordinary double point if and only if $s(a)$ is a square (resp.\ a
    nonsquare) in $\FF^\times$. 
  \item
    The normalization $\Yb_0$ of $\Yb$ is given by the equation
    \[
        \vb^2 = s(\xb).
    \]
    The map $\Yb_0\to\Yb$ is determined by $\ub=r\vb$.
  \end{enumerate}
\end{prop}
 
\proof The proof is similar to but easier than the proof of Proposition
\ref{prop:redp=2}, and is therefore omitted.  \Endproof

\begin{cor} \label{cor:oddp}
  Assume that $\Yb$ is semistable. 
  \begin{enumerate}
  \item
    There is a bijection between the set $\Sc$ of singular points of $\Yb$ and
    the irreducible factors of the polynomial $r$. 
  \item
    For $\xi=(a,b)\in\Sc$ we have $\varepsilon_\xi=1$ (resp.\
    $\varepsilon_\xi=-1$) if and only if $s(a)$ is
    a square (resp.\ a nonsquare) in $\FF_\p^\times$.
  \item
    The conductor exponent is
    \[
       f_\p = \deg (r).
    \]
  \end{enumerate}
\end{cor}

\subsection{} 

In this section we summarize the results obtained so far and describe
the resulting algorithm for computing the $L$-function of the curve $Y$.
We also make some comments on the implementation and running time. For
simplicity we assume from now on that $K=\QQ$.

We are interested in computing a certain $L$-series given as an Euler product,
\[
     L(Y/K,s) = \sum_{n\geq 1} \frac{a_n}{n^s} = \prod_p L_p(Y,s).
\]
More specifically, we want to give evidence for the conjectured functional
equation (Equation (\ref{eqn:funceq})). We use the \texttt{Dokchitser} package
in the free computer algebra software \textsl{sage}, based on Tim Dokchitser's
paper \cite{Dokchitser04}. To verify the functional equation in this package,
we need to know the conductor of the $L$-function,
\[
       N = \prod_p p^{f_p},
\]
and the coefficients $a_n$ need to be calculated for all $n\leq M$ up to a
certain heuristic constant $M$, depending on $N$ and $g_Y$. The constant $M$
can be computed using the \texttt{Dokchitser} package. 
Due to the multiplicativity relation of the $a_n$, we only have to compute the
coefficients $a_{p^k}$ for prime powers $p^k\leq M$ via point counting. If one
uses naive point counting methods, the calculation of each $a_{p^k}$ has a
complexity of about $\mathcal O(p^k)$. By the prime
number theorem, we get a complexity of $\mathcal O(M/2\cdot M/\log(M))$ for
each $L$-series.
For fixed genus, $M$ is proportional to the square root of the conductor $N$
of the curve (cf. \cite{Dokchitser04}). So the
complexity of checking the functional equation is bounded by $\mathcal
O(N/\log N)$.

Finding examples of suitable curves $Y$ and checking the functional equation
of $L(Y,s)$ can be performed as follows.\label{sec:algo} We fix an integer
$g_Y\geq 2$.

\begin{itemize}
\item[1.]  Choose polynomials $g,h\in\ZZ[x]$ with $h\not \equiv
  0\pmod{2}$, $g$ monic, $\deg(g)=2g_Y+1$, $\deg(h)\leq g_Y$.  Consider the
  reductions of $g,h$ modulo 2. If 
  \[
      \gcd(\bar h,\bar h',\bar g')\neq 1,
  \] 
  the model $\Y$ from \S~ \ref{ssred2} is not semistable at $p=2$ (Corollary
  \ref{cor:redp=2}). In this case we dismiss our choice of $g$ and $h$ and
  start over again. Otherwise, compute the polynomial
  \[
      r:=\gcd(\hb,(\hb')^2\gb+(\gb')^2).
  \]
\item[2.]  Calculate the discriminant $\Delta\in\ZZ$ of the polynomial
  $f:=4g+h^2$, and define $S'$ as the set of prime factors of $\Delta$,
  ignoring the prime factor $2$. Check for all $p\in S'$ whether
  \[
      \gcd(\fb,\fb',\fb'')=1.
  \]
  If the test fails for one $p\in S'$ then we cannot guarantee that $Y$ has
  semistable reduction (Lemma \ref{lem:podd}). If this happens we dismiss
  our example and go back to the beginning. If $\deg r>1$ then we set
  $S:=S'\cup\{2\}$, otherwise set $S:=S'$.

  Now we know that $Y$ has bad semistable reduction at all primes $p\in S$ and
  good reduction everywhere else. 
\item[3.] For all \textsl{bad primes} $p\in S$, we do the following.
  \begin{itemize}
  \item[3a.] For $p=2$, decompose the polynomial $r\in\FF_2[\xb]$ into
    irreducible factors,
    \[
        r=\prod_i r_i.
    \]
    Each factor $r_i$ corresponds to a singularity $\xi_i\in\Yb$ with
    $\deg(\xi_i)=\deg(r_i)$. Check for all $i$ whether $\xi_i$ is split or not
    (Corollary \ref{cor:redp=2.2} (ii)) and set $\varepsilon_i\in\{\pm 1\}$
    accordingly.

    Now calculate the numerator $P(\Yb_0,T)$ of the zeta function of the
    normalization $\Yb_0$ of $\Yb$, using Equation \eqref{eqn:redp=2.3} and
    point counting (Remark \ref{rem:pointcount} and
  Remark \ref{rem:pointcount2}). The local $L$-factor at
    $p=2$ is
    \[
       L_2(\Yb,T)= P(\Yb_0,T)^{-1}\prod_i (1-\varepsilon_iT^{d_i})^{-1},
    \]
    see Proposition \ref{prop:badfactor}. Also, $f_2=\deg(r)$ (Corollary
    \ref{cor:redp=2.2} (iii)).

  \item[3b.] For $p\in S-\{2\}$, compute the decomposition
    \[
          \fb = r^2 s.
    \]
    Split $r=\prod _i r_i\in \FF_p[x]$ into irreducible factors. Set
    $d_i:=\deg(r_i)$ and $\varepsilon_i:=\pm 1$, according to Corollary
    \ref{cor:oddp} (ii).

    Calculate the numerator $P(\Yb_0,T)$ of the zeta function of the
    normalization $\Yb_0$ of $\Yb$, using the equation from Proposition
    \ref{prop:podd} (iii). As in 3a, the local $L$-factor is
    \[
       L_p(\Yb,T)= P(\Yb_0,T)^{-1}\prod_i (1-\varepsilon_iT^{d_i})^{-1}.
    \]
    Set $f_p:=\deg(r)$ (Corollary \ref{cor:oddp} (iii)). 
  \end{itemize}
\item[4.]
  Compute the conductor
  \[
        N:=\prod_{p\in S} p^{f_p}
  \]
  and the constant $M$.
\item[5.] Calculate the local $L$-factor $L_p(Y,s)$ for all good
  primes $p\notin S$, $p\leq M$ via point counting (Remarks
  \ref{rem:pointcount} and \ref{rem:pointcount2}). 
\item[6.] Compute the truncated $L$-series
  \[
       L(Y,s)'=\displaystyle\sum\limits_{n=1}^M\frac{a_n}{n^{-s}}
  \]
  from the Euler factors $L_p(Y,s)$, $p\leq M$. Check the functional equation
  using the \texttt{Dokchitser} package for the root number $1$. If this
  fails, repeat with root number $-1$. 
\end{itemize}

\begin{rem}\label{rem:trunc}
  The algorithm described above can be slightly improved as follows.  Observe
  that we only need coefficients $a_n$ of the $L$-series with $n\leq M$. Thus
  we can use the bound $M$ (which only depends on the conductor exponents
  $f_p$ for $p\in S$) to truncate the polynomial $P(\Yb_0,T)$, resp. the local
  $L$-factor $L_p$, in order to avoid counting points over fields with more
  than $M$ elements. This affects the computation of the polynomials
  $P(\Yb_0,T)$ in Step 3a and 3b.  For an example, we refer to Example
  \ref{exa:trunc} in \S~\ref{sec:exa}.
\end{rem}

In the above algorithm, the time needed to compute the set of bad primes, the
conductor $N$ and the constant $M$ is insignificant compared to the time
needed for the point counting. The numerical verification of the
functional equation is not expensive either. Therefore, the running time of
our algorithm for an individual curve $Y$ is indeed bounded by $\mathcal
O(N/\log N)$, as explained above, with almost all the running time spent on
point counting. For the class of hyperelliptic curves considered in this
section, examples with conductor up to $10^{10}$ can be computed within a
reasonable time. In the largest example that we computed, the conductor is
$N=7 \cdot 11 \cdot 13 \cdot 89 \cdot 431 \cdot 857 \approx
3\cdot10^{10}$. Using more sophisticated point-counting methods as e.g.\ in
\cite{KedlayaSutherland} would probably allow the computation of significantly
larger examples.

The running time of our example is essentially determined by the
conductor. Although the constant $M$ depends on $N$ and the genus of $Y$, its
dependence on $g_Y$ is insignificant within the range of genera that we
consider (an asymptotic estimate for $M$ can be obtained from
\cite{Dokchitser04}, \S 4, in particular Equation (4-2)). This is an advantage
of our approach, as opposed to for example that of Booker (\cite{booker},
\S~2.3.2). However, the discriminant of the polynomial $f=4g+h^2$ determines
the odd prime factors of the conductor. Thus with growing degree of $g$, and
therefore with growing genus $g_Y=(\deg (g)-1)/2$, it gets more and more
difficult to find examples of curves $Y$ with conductor of reasonable size. So
far, we managed to find examples that fall within this range for all $g_Y\leq
6$.

We have verified the functional equation for several hundreds of
examples.  Obviously, one can easily construct a lot more examples,
especially for small genus. On our homepage, we provide a selection of
examples -- each with slightly different parameters -- where the
functional equation has been verified. The data can be found on {\tt
  https://www.uni-ulm.de/index.php?id=64504}.

\subsection{}\label{sec:exa}

In this section we give a few explicit examples, in detail.  All given
examples have been checked to fulfill the functional equation. Note that the
chosen examples do not necessarily have the smallest possible conductor for
the given genus -- it is merely a selection of `typical' examples.

\begin{exa}
  The polynomials 
  \[
        g = x^5 - 3x^4 - 3x^3 - 3x^2 - 3x - 1, \quad 
        h=x^2 + 3x + 1
  \] 
  define a genus-two curve $Y/\QQ$. We find five {bad primes}:
  $2,3,7,101,163$. The $L$-factors corresponding to these primes are (we write
  $T$ instead of $p^{-s}$):
  \begin{align*}
    L_2^{-1}&=1+T^2,\\
    L_3^{-1}&=(1+T)  (3T^2 - T + 1),\\
    L_7^{-1}&=(1-T) (7T^2 + 3T + 1),\\
    L_{101}^{-1}&=(1+T)  (101T^2 + 3T + 1),\\
    L_{163}^{-1}&=(1-T) (163T^2 + 11T + 1).
  \end{align*}
  The conductor is $N=2^2 \cdot 3 \cdot 7 \cdot 101 \cdot 163\approx 10^6$.
  We briefly review the computations for $p=2,3$. For
  $p=2$ we look at the curve
  \[
        \Yb/\FF_2:\; \yb^2+(1+\xb+\xb^2)\yb = 1+\xb+\xb^2+\xb^3+\xb^4+\xb^5.
  \]
  Since $\hb'=(1+\xb+\xb^2)'=1$, $\Yb$ is semistable. The singular locus is
  determined by the polynomial
  \[
     r:=\gcd(\hb,(\hb')^2\gb+(\gb')^2)=\hb=1+\xb+\xb^2.
  \]
  Hence there is a unique ordinary double point $\xi=(a,b)$ of degree $2$,
  where $a$ is a solution to $a^2+a+1=0$. Substituting $\yb=\hb\tilde{y}$ into
  the equation for $\Yb$ and dividing by $\hb^2$ we obtain the equation for
  its normalization:
  \[
       \tilde{y}^2+\tilde{y} = \gb/\hb^2=1+\xb.
  \]
  This is a curve of genus zero, so it doesn't contribute to the local
  $L$-factor. However, we can see that the inverse image $\pi^{-1}(\xi)$ of
  the singular point $\xi$ corresponds to the solutions to the equation
  \[
        \tilde{y}^2+\tilde{y}=1+a
  \]
  in $\FF_2(a)=\FF_4$. Clearly, this equation is irreducible and so $\xi$ is a
  non-split ordinary double point. Therefore,
  \[
          L_2(Y,s) = \frac{1}{1+2^{-2s}}.
  \] 

  For $p=3$ we are looking at the curve
  \[
     \Yb/\FF_3:\; \bar{u}^2=\xb^2(2+\xb^2+\xb^3).
  \]
  This is a semistable curve with one $\FF_3$-rational ordinary double point
  $\xi=(0,0)$. Substituting $\ub=\xb\vb$ and dividing by $\xb^2$ gives an
  equation for the normalization of $\Yb$,
  \[
      \Yb_0:\; \vb^2=2+\xb^2+\xb^3,
  \]
  a smooth curve of genus $1$ over $\FF_3$. There are exactly three rational
  points, $\abs{\Yb_0(\FF_3)}=3$. It follows that the numerator of the zeta
  function is $P(\Yb_0,T)=1-T+3T^2$. Also, the fiber $\pi^{-1}(\xi)$ is given
  by the equation
  \[
       \vb^2=2,
  \]
  which is irreducible over $\FF_3$. It follows that $\xi$ is a non-split
  double point and that
  \[
       L_3(Y,s) = \frac{1}{(1+3^{-s})(1-3^{-s}+3^{1-2s})}.
  \]
  The computation of the $L$-factors for $p=7,101,163$ is similar.

  The numerical verification of the functional equation was successful, with
  root number 1.
\end{exa}

\begin{exa}
  The polynomials $g= x^7 + x^6 + 2x^5 + 2x^4 + 2x^3 - 1$ and $h= -x^3
  + x^2 + x + 2$ define a genus-three curve. We find four {bad
    primes}: $2,3,11,37$. The $L$-factors corresponding to these primes
  are
  \begin{align*}
    L_2^{-1}&=(1-T)  (1+T)  (2T^2 - T + 1),\\
    L_3^{-1}&=(1+T) (1-T^2), \\
    %\quad,\quad \bar r=(x + 1)  (x^2 + 1)\\
    L_{11}^{-1}&=( 1+T)^2 (11 T^2 - 4 T + 1), \\
    L_{37}^{-1}&=(1-T)(37^2 T^4 + 148 T^3 + 14 T^2 + 4 T + 1), 
  \end{align*}
  the conductor is $N=2^2\cdot3^3\cdot 11^2\cdot 37 \approx 10^5$, and the
  root number is $1$.
\end{exa}

\begin{exa}The polynomials $g= x^9 - 2x^8 + x^7 - 2x^4 + 2x^3 + 2x^2 + x$ and
  $h= -2x^4 + x^3 - 2x^2 - x - 1$ define a genus-four curve. We find four
  {bad primes}: $3,7,31,53$. The $L$-factors corresponding to these
  primes are
  \begin{align*} L_3^{-1}&=(1-T) (1+T) (9T^4 + 6T^3 + 4T^2 + 2T + 1), \\
    L_7^{-1}&=(1 + T^3) (7T^2 + 3T + 1), \\
    L_{31}^{-1}&=(1 - T) (31^3T^6 + 1581T^4 + 36T^3 + 51T^2 + 1), \\
    L_{53}^{-1}&=(1+T) (53^3T^6 + 8427T^5
    +1537T^4 + 670T^3 + 29T^2 + 3T + 1), 
  \end{align*}
  the conductor is $N=3^2 \cdot 7^3 \cdot 31 \cdot 53 \approx 10^8$, and the
  root number is $1$.
\end{exa}

\begin{exa} 
  The polynomials $g= x^{11} + 3x^4 + 2x^3 - 3x^2 - 2x$ and $h=
  -3x^3\label{exa:trunc} + x^2 + 3x + 1$ define a genus-five curve. We find
  four {bad primes}: $7,227,1277,1609$. The truncated $L$-factors
  corresponding to these primes are
  \begin{align*} 
    L_7^{-1}&=(1-T)( 7^4T^8 - 588T^6 + 134T^4 - 12T^2 + 1), \\
    L_{227}^{-1}&=(1+T) (\ldots + 200T^2 + 13T + 1), \\
    L_{1277}^{-1}&=(1+T)  (\ldots-35T+1),\\
    L_{1609}^{-1}&=(1+T) (\ldots-26T+1).
  \end{align*}
  The conductor is $N= 7 \cdot 227 \cdot 1277 \cdot 1609 \approx
  10^{9}$, and the root number is $1$.  We truncated the last three
  $L$-factors to save computation time, since the bound in this example
  is $M=1112661<\{227^3,1277^2,1609^2\}$. Hence no further information
  is needed to verify the functional equation (\S~ \ref{sec:algo}).
\end{exa}

\begin{exa} 
  The polynomials $g= x^{13} + x^{12} + x^{11} + x^{10} + x^9 + x^8 - x^7 - x^6 + x^5 + x^4 + x^3$ and $h=
  x^6 + x^5 - x^3 + x^2 + x + 1$ define a genus-six curve. We find
  six {bad primes}: $7, 11 , 13 , 89 , 431 , 857$. The truncated $L$-factors
  corresponding to these primes are
  \begin{align*} 
    L_7^{-1}&=(1+T)(7^5T^{10} + 1372T^8 + 1127T^7 + 112T^6 + 122T^5 +
    16T^4 + 23T^3 +\\
&\quad 4T^2 + 1 ), \\ L_{11}^{-1}&=(1+T) (11^5T^{10} -
    43923T^9 - 10648T^8 + 2662T^7 + 781T^6 - 390T^5 + \\
&\quad 71T^4 + 22T^3 -
    8T^2 - 3T + 1) , \\ L_{13}^{-1}&=(1+T) (13^5T^{10} - 114244T^9 +
    43940T^8 - 10140T^7 + 1040T^6+\\
&\quad - 342T^5 + 80T^4 - 60T^3 + 20T^2 -
    4T + 1) ,\\ L_{89}^{-1}&=(1-T) (\ldots+320T^3 + 43T^2 - 5T +
    1),\\ L_{431}^{-1}&=(1+T) (\ldots+859T^2 + 31T + 1
    ),\\ L_{857}^{-1}&=(1-T) (\ldots+1448T^2 - 41T + 1 ) \end{align*}
  The conductor is $N= 7 \cdot 11 \cdot 13 \cdot 89 \cdot 431 \cdot
  857 \approx 3\cdot10^{10}$, and the root number is $1$.
  As in the previous example, we truncated the last three $L$-factors,
  for the bound is $M=2549728<\{ 89^4, 431^3, 857^3\}$.
\end{exa}

\section{More examples}\label{sec:more_examples}

Combining the results of \cite{superell}, \cite{ArzdorfWewers},
\cite{ArzdorfDiss} and \cite{superp} it is in principle possible to
compute the semistable reduction of any superelliptic curve and therefore the
local $L$-factors and the conductor exponents, at all primes. In the
previous section we have chosen a class of examples where this was
particularly easy. In this section we treat a small sample of examples which do
not fall within this class, but where we were nevertheless able to compute
$L_p(Y,s)$ and $f_p$ for all $p$ and to verify the functional equation.

\subsection{}

The polynomials 
\[
     g= x^7 - 2x^6 - 2x^4 + x^3 + 3x^2 + x,\quad 
     h= 3x^3 + 3x^2 + 2x + 1
\]
substituted in (\ref{eqn:p=2})
define a genus-three hyperelliptic curve $Y/\QQ$. One checks that $Y$ has good
reduction at $p=2$ (see Step 1 and 2 in the algorithm from \S~
\ref{sec:algo}). We ignore the prime $p=2$ from now on and use the equation
\begin{equation} \label{eqn:exa1.1}
    Y/\QQ:\; y^2=f(x)=4x^7 + x^6 + 18x^5 + 13x^4 + 22x^3 + 22x^2 + 8x + 1,
\end{equation}
with $f:=h^2+4g$, to describe $Y$. (As in \S~ \ref{sec:ssred} this means that
the open affine subset $Y-\{\infty\}\subset Y$ is the plane affine curve given
by the above equation.) The discriminant of $f$ is $\Delta=-2^{12}\cdot3\cdot
5^3\cdot 13^2\cdot 97$. Therefore, there are four bad primes:
$p=3,5,13,97$. For the primes $p=3,13,97$ the condition
$\gcd(\fb,\fb',\fb'')=1$ holds and therefore $Y$ has semistable reduction at
$p$. The local $L$-factors and the conductor exponent can be computed as
before. We obtain
\begin{align*}
  L_3^{-1}&=(1-T) (3^2T^4 + 3T^3 + T + 1), \quad f_3=1,\\
  L_{13}^{-1}&=(1+T^2)  (13T^2 + 5T + 1),\quad f_{13}=2,\\
  L_{97}^{-1}&=(1+T) (97^2T^4 + 582T^3 + 78T^2 + 6T + 1), \quad f_{97}=1.
\end{align*}
However, for $p=5$ the special fiber of the naive model $\Y^{\rm naive}$ of $Y$
over $\ZZ_{(5)}$ (obtained by reducing equation \eqref{eqn:exa1.1} modulo $5$)
is the curve
\begin{equation} \label{eqn:exa1.2}
    \Yb^{\rm naive}/\FF_5:\; \yb^2=\fb=4(\xb+1)^4(\xb^3+\xb+4).
\end{equation}
We see that $\Yb^{\rm naive}$ has a unique $\FF_5$-rational singularity
$(\xb,\yb)=(4,0)$ which is not an ordinary double point, and is smooth
everywhere else. In particular, $\Yb^{\rm naive}$ is not semistable, and the
methods from \S ~\ref{sec:ssred} are not directly applicable. 

Nevertheless, using the results of \cite{superell} 
we can easily compute the semistable reduction of $Y$ at
$p=5$.  
 We are dealing with
a local problem and may therefore consider $Y$ as a curve over the
$5$-adic numbers $\QQ_5$. Let $L:=\QQ_5[\pi]$ be the extension of
degree $4$ with $\pi^4=5$. Clearly $L/\QQ_5$ is a Galois extension,
which is totally and tamely ramified. The Galois group of $L/\QQ_5$ is cyclic,
generated by the element $\sigma$ determined by
\[
       \sigma(\pi)=\zeta_4\pi.
\]
Here $\zeta_4\in\ZZ_5$ is the $4$th root of unity with $\zeta_4\equiv
2\pmod{5}$. Let $\p=(\pi)\lhd\OO_L$ denote the unique prime ideal. Note that
the residue field is $\FF_\p=\FF_5$.

\begin{lem} \label{lem:exa1}
\begin{enumerate}
\item
  The curve $Y_L=Y\otimes_{\QQ_5} L$ has semistable reduction at $\p$.
\item Let $\Y/\OO_L$ denote the minimal semistable model of $Y_L$ and
  $\Yb$ its special fiber. The curve $\Yb$ is the union of two smooth,
  absolutely irreducible curves over $\FF_5$. The first component
  $\Yb_1$ has an affine open subset which is given by the equation
  \begin{equation} \label{eqn:exa1.3}
      \yb_1^2=4(\xb^3+\xb+4),
  \end{equation}
  the second component $\Yb_2$ has an affine open subset given by
  \begin{equation} \label{eqn:exa1.4}
      \yb_2^2=3\xb_2^4+2.
  \end{equation}
\item
  The components $\Yb_1,\Yb_2$ intersect in a unique split ordinary double
  point $\xi$ of degree $2$. As a point on $\Yb_1$, we have $\xi=(4,b)$, where
  $b\in\FF_5(\xi)=\FF_{5^2}$ is a solution to $b^2=3$.  
\item
  The generator $\sigma$ of the Galois group
  $\Gal(L/\QQ_5)=\gen{\sigma}$ acts trivially on $\Yb_1$ and on
  $\Yb_2$ via the automorphism of order $4$
  \[
        \xb_2\mapsto 3\xb_2, \quad \yb_2\mapsto -\yb_2.
  \]
\end{enumerate}
\end{lem}

\proof
One simply follows the recipe in \cite{superell}, \S~4. The equations for
$\Yb_i$, $i=1,2$, are obtained as follows. For $\Yb_1$, we substitute
$y=(x+1)y_1$ in \eqref{eqn:exa1.1}, divide by $(x+1)^4$ and reduce modulo
$\p$. Using \eqref{eqn:exa1.2} we see that we obtain the equation
\ref{eqn:exa1.3}. For $\Yb_2$, we substitute $x=4+\pi x_1$ and $y=\pi^2y_2$,
divide by $5$ and reduce modulo $\p$. A short computation yields 
\eqref{eqn:exa1.4}.  Statements (iii) and (iv) are straightforward.
\Endproof

We note that by construction the semistable model $\Y$ dominates the naive
model $\Y^{\rm naive}$, or in other words, there is a modification
$\Y\to\Y^{\rm naive}\otimes_{\ZZ_5}\OO_L$. The resulting map
$\Yb\to\Yb^{\rm naive}$ may be visualized as in Figure \ref{fig:Yb_4_1}.

\begin{figure}[h!]\centering 
\begin{tikzpicture}[scale = 0.12]
 \coordinate (A) at (-25,2);
 \coordinate (G) at (-1.8,5);
 \coordinate (G0) at (-2,12);
 \coordinate (H) at (1,14.8);

 \coordinate (I0) at (4,12);
 \coordinate (I) at (3.8,5);
 \coordinate (Z) at (30,2);

 \coordinate (C) at (-25,25);
 \coordinate (D) at (-10,23);
 \coordinate (E) at (16,27);
 \coordinate (F) at (30,25);
 \coordinate (P) at (-4,20);
 \coordinate (Q) at (-2,32);
 \coordinate (R) at (4,32);
 \coordinate (S) at (6,20);

%singular curve
\draw [black] plot [smooth, tension=1.0] coordinates {   (A)(G)(G0)(H)} ;
\draw [black] plot [smooth, tension=1.0] coordinates {   (Z)(I)(I0)(H)} ;
\node [right] at (Z) {$\bar Y^{\rm naive}$};

% normalization part 1
\draw [black] plot [smooth, tension=1.0] coordinates {   (C)(D)(E)(F)} ;
\node [right] at (F) {$\bar Y_1$};

% normalization part 2
\draw [black] plot [smooth, tension=1.0] coordinates {   (P)(Q)(R)(S)} ;
\node [right] at (R) {$\bar Y_2$};
% Hilfslinien
%\draw [gray!77] plot coordinates {   (A)(G)(G0)(H)} ;
%\draw [gray!77] plot coordinates {   (Z)(I)(I0)(H)} ;
\end{tikzpicture}
{\caption{$\Yb\to\Yb^{\rm naive}$}\label{fig:Yb_4_1}}
\end{figure}

We see that $\Yb_2$ is contracted to the singular point on $\Yb^{\rm naive}$ and
that $\Yb_1$ can be identified with the normalization of $\Yb^{\rm naive}$. 

\begin{cor} \label{cor:exa1}
  The local $L$-factor and the conductor exponent of the curve $Y$ at $p=5$
  are 
  \[
     L_5^{-1}=(1+T) (1+3T+5T^2),\quad f_5=3.
  \]
\end{cor}

\proof Let $\Zb:=\Yb/\Gal(L/\QQ_5)$ be the quotient of $\Yb$ under the action
of the Galois group of the extension $L/\QQ_5$. In the terminology of
\cite{superell}, $\Zb$ is the {\em inertial reduction} of $Y$ at $p=5$. It
follows from Lemma \ref{lem:exa1} that $\Zb$ is a semistable curve over
$\FF_5$, consisting of two irreducible components $\Zb_1,\Zb_2$ which
intersect in a unique non-split ordinary double point of degree $2$. The curve
$\Zb_1$ may be identified with $\Yb_1$ and the curve $\Zb_2$ with the quotient
$\Yb_2/\gen{\sigma}$. One sees immediately from \eqref{eqn:exa1.4} and Lemma
\ref{lem:exa1} (iv) that $\Zb_2$ has genus zero. 

By \cite{superell}, Corollary 2.5, the local $L$-factor is
\[
     L_5(Y,s) = P(\Zb,5^{-s})^{-1},
\]
where $P(\Zb,T)$ is the numerator of the zeta function of $\Zb$. From the
above description of $\Zb$ we see that 
\[
     P(\Zb,T) = (1+T)(1+3T+5T^2).
\]
The second factor is the numerator of the zeta function of the genus-one curve
$\Zb_1$, which is given by \eqref{eqn:exa1.3}, and the first factor comes from
the action of $\Gamma_{\FF_5}$ on $H^1(\Delta_{\Zb})$, as in the proof of
Proposition \ref{prop:badfactor}. Finally, we use \cite{superell}, Corollary
2.6, to conclude that
\[
    f_5= 2g_Y-\dim H^1_\et(\Zb_k,\QQ_\ell) = 6 - 3 = 3.
\]
\Endproof

We have computed the local $L$-factors and conductor exponents at all
bad primes. We can now continue with Step 4 of the algorithm from
\S~\ref{sec:algo}. The conductor of the $L$-function is $N= 3 \cdot
5^3 \cdot 13^2\cdot 97 \approx 10^8$. We set $M:=55956$, compute the
local $L$-factors for all good primes $p\leq M$ and the truncated
$L$-series $L(Y,s)'$. Feeding these data into Dokchitser's algorithm,
we have checked that the $L$-function of $Y$ verifies the expected
functional equation with root number $-1$.

\subsection{}

We now treat an example of a hyperelliptic curve which does not have semistable
reduction at $p=2$. In this case, the methods of \cite{superell} to compute
semistable reduction do not apply. 

The polynomials 
\[
  g= x^9 - x^8 + x^7 + x^5 + x^3,\quad
  h= -x^4 + 1
\]
define a hyperelliptic curve of genus four. The discriminant of $f:=h^2+4g$ is
$\Delta=-2^{32}\cdot 317$. So $p=317$ is the only odd prime where $Y$ has bad
reduction. Running through Step 2 and 3b of the algorithm from
\S~\ref{sec:algo} we see that $Y$ has semistable reduction at $p=317$ and the
local $L$-factor and the conductor exponent are
\[
    L_{317}^{-1}=(1+T) (1-32T+991T^2+\ldots), \quad f_{317}=1.
\]
It will follow from the calculation of the conductor $N$ below that
this is indeed the correct truncation (Remark \ref{rem:trunc}).

Unfortunately, the condition $\gcd(\hb,\hb',\gb')=1$ from \S~\ref{sec:algo},
Step 1, is not satisfied. The naive model $\Y^{\rm naive}$ of $Y$ over $\ZZ_{(2)}$,
given by \eqref{eqn:p=2}, has special fiber
\begin{equation} \label{eqn:exa2.1}
      \Yb^{\rm naive}/\FF_2:\; \yb^2+(\xb+1)^4\yb = \gb(\xb)=
          \xb^9+\xb^8+\xb^7+\xb^5+\xb^3.
\end{equation}
One sees at once that $\Yb^{\rm naive}$ has a non-ordinary singularity at
$(\xb,\yb)=(1,1)$. Substituting 
\[
     \yb=(\xb+1)^3\yb_0+\xb^4+\xb^2+1
\]
into \eqref{eqn:exa2.1} and dividing by $(\xb+1)^6$ we obtain an equation for
the normalization $\Yb_0'$ of $\Yb^{\rm naive}$:
\begin{equation} \label{eqn:exa2.2}
  \Yb_0'/\FF_2:\qquad \yb_0^2+(\xb+1)\yb_0 = \xb^2(\xb+1).
\end{equation}
We see that $\Yb_0'$ is a smooth curve of genus one. The numerator of
its zeta function is
\begin{equation} \label{eqn:exa2.3}
     P(\Yb_0',T) = 1+T+2T^2.
\end{equation}

The computation of the semistable reduction of $Y$ at $p=2$ is rather
challenging. We only state the result (Lemma \ref{lem:exa4.2} below). A
detailed proof will be given elsewhere. 

Let us work over the field $\QQ_2$ of $2$-adic numbers. Using the methods of
\cite{superp} we produce the following polynomial:
\begin{multline}
  \Delta=x^{12} + 20 x^{11} + 154 x^{10} + 664 x^9 + 1873 x^8 + 3808 x^7
       + 5980 x^6 \\+ 7560 x^5 + 7799 x^4 + 6508 x^3 + 4290 x^2 + 2224 x + 887
       \in \ZZ_2[x].
\end{multline}
One  checks that $\Delta$ is irreducible over $\QQ_2$. Let $L/\QQ_2$ be
the splitting field of $\Delta$, $\Gamma=\Gal(L/\QQ_2)$ the Galois group and
$K_i:=L^{\Gamma_i}$ the fixed field of the $i$th ramification group, for
$i\geq 0$. One also checks that $K_0/\QQ_2$ has degree $2$ and that $K_1/K_0$ has
degree $9$. So $K_0/\QQ_2$ is the unique unramified extension of degree $2$,
and $\Gamma_0/\Gamma_1$ is a cyclic group of order $9$. Unfortunately, we do
not know the exact size and structure of the wild inertia group
$\Gamma_1$. Nevertheless, we can prove the following.

\begin{lem} \label{lem:exa4.2}
\begin{enumerate}
\item
  The curve $Y_L$ has semistable reduction.
\item Let $\Y$ be the minimal semistable model of $Y$ over $\OO_L$ and $\Yb$
  the special fiber of $\Y$. Then $\Yb$ consists of five irreducible
  components $\Yb_0,\ldots,\Yb_4$ over the residue field $\FF_\p=\FF_4$ of
  $L$. Here $\Yb_0$ may be identified with the pullback to $\FF_\p$ of the
  curve $\Yb_0'$, the normalization of $\Yb^{\rm naiv}$. The components
  $\Yb_1,\Yb_2,\Yb_3$ are smooth curves of genus one over $\FF_\p$, given by
  equations
  \[
       \Yb_i/\FF_p:\;\yb_i^2+\yb_i = \xb_i^3, \qquad i=1,2,3.
  \]
  The component $\Yb_4$ is a projective line and intersects each of the other
  four components in a unique point. The genus one components
  $\Yb_0,\ldots,\Yb_3$ do not intersect (Figure \ref{fig:Yb_4_2}).
\item The inertia group $\Gamma_0$ fixes $\Yb_0$ and $\Yb_4$ and permutes the
  components $\Yb_1,\Yb_2,\Yb_3$ transitively. The wild inertia group
  $\Gamma_1$ fixes every component.
\item Let $\Gamma_0'\subset \Gamma_0$ be the stabilizer of the component
  $\Yb_1$, $H\subset \Gamma_0'$ the kernel of the map
  $\Gamma_0'\to\Aut(\Yb_1)$, and $\tilde{\Gamma}_0=\Gamma_0'/H$ the
  quotient. Then $\tilde{\Gamma}_0$ is cyclic of order $6$. Its unique element
  of order two acts on $\Yb_1$ via the automorphism
  \[ 
      \xb_1\mapsto \xb_1, \quad \yb_1\mapsto \yb_1+1.
  \]
  Moreover, the filtration of higher ramification groups on $\tilde{\Gamma}_0$
  has the form
  \[
      \tilde{\Gamma}_0 \supsetneq \tilde{\Gamma}_1=\ldots=\tilde{\Gamma}_{15}
                        \supsetneq \tilde{\Gamma}_{16}=1.
  \]
\end{enumerate}
\end{lem}

\begin{figure}[h!]\centering 
\begin{tikzpicture}[scale = 0.12]
 \coordinate (A) at (-15,-13);
 \coordinate (A1) at (-26,-12);
 \coordinate (B) at (2,-10);
 \coordinate (B1) at (11,-9);
 \coordinate (C) at (40,-10);
 \coordinate (C1) at (45,-8);
 \coordinate (D) at (-20,-25);
 \coordinate (D1) at (-10,-25);
 \coordinate (D2) at (-5,-25);
 \coordinate (E0) at (2,-9.5);
 \coordinate (E00) at (2,-12);
 \coordinate (E) at (0,-20);
 \coordinate (E1) at (2,-5);
 \coordinate (E2) at (4,-20);
 \coordinate (E3) at (12,-23);
 \coordinate (E4) at (25,-25);
 \coordinate (F) at (40,-24);
 \coordinate (F1) at (43,-24);
 \coordinate (F2) at (53,10);
 \coordinate (G0) at (-4,-2);
 \coordinate (G) at (2,-3);
 \coordinate (G1) at (15,-1);
 \coordinate (G2) at (25,-3);
 \coordinate (H0) at (-4,-2);
 \coordinate (H) at (2,-3);
 \coordinate (H1) at (15,-1);
 \coordinate (H2) at (25,-3);
 \coordinate (H3) at (30,-1);
 \coordinate (I0) at (-4,3);
 \coordinate (I) at (2,2);
 \coordinate (I1) at (15,4);
 \coordinate (I2) at (25,2);
 \coordinate (I3) at (30,4);
 \coordinate (J0) at (-4,8);
 \coordinate (J) at (2,7);
 \coordinate (J1) at (15,9);
 \coordinate (J2) at (25,7);
 \coordinate (J3) at (30,9);
 \coordinate (K) at (2,10);

%singular curve left and right part
\draw [black] plot [smooth, tension=1.0] coordinates {   (D)(E)(E1)} ;
\draw [black] plot [smooth, tension=1.0] coordinates {   (E1)(E2)(E4)(F)} ;
\node [right] at (F1) {$\bar Y^{\rm naive}$};

% delete the loop
 \filldraw [gray!00]      (E0) circle (226pt) ;

%  horizontal curve
\draw [black] plot [smooth, tension=1.0] coordinates {   (A1)(A)(B1)(C)} ;
\node [below] at (C1) {$\Yb_0$};
%\node [below] at (F2) {$\bar Y$};

% vertical line
\draw [black] plot [smooth, tension=1.0] coordinates {   (E00)(K)} ;

% short horizontal curves
%\draw [black] plot [smooth, tension=1.0] coordinates {   (G0)(G)(G1)(G2)} ;

\draw [black] plot [smooth, tension=1.0] coordinates {   (H0)(H)(H1)(H2)} ;
\node [below] at (H3) {$\Yb_1$};

\draw [black] plot [smooth, tension=1.0] coordinates {   (I0)(I)(I1)(I2)} ;
\node [below] at (I3) {$\Yb_2$};

\draw [black] plot [smooth, tension=1.0] coordinates {   (J0)(J)(J1)(J2)} ;
\node [below] at (J3) {$\Yb_3$};

\end{tikzpicture}
{\caption{$\Yb\to\Yb^{\rm naive}$}\label{fig:Yb_4_2}}
\end{figure}

\begin{cor} \label{cor:exa2}
  The local $L$-factor and the conductor exponent of $Y$ at $p=2$ are
  \[
      L_2^{-1}=1+T+2T^2, \qquad f_2=16.
  \]
\end{cor} 

\proof Let $\Zb=\Yb/\Gamma$ be the quotient curve. It follows directly from
Lemma \ref{lem:exa4.2} that $\Zb$ is a semistable curve over $\FF_2$,
consisting of three irreducible components (corresponding to the three orbits
of the action of $\Gamma_1$ on the set of irreducible components of
$\Yb$). The first component is the genus one curve
$\Zb_0:=\Yb_0/\Gamma\cong\Yb_0'$, given by
\eqref{eqn:exa2.2}. The other two components have genus zero. Moreover, the
component graph of $\Zb$ is a tree. It follows that the zeta function of $\Zb$
is the same as the zeta function of $\Yb_0'$, and hence
\[
        P(\Zb,T)=P(\Yb_0',T)=1+T+2T^2,
\]
by \eqref{eqn:exa2.3}. The claim $L_2^{-1}=1+T+2T^2$ follows now from
\cite{superell}, Corollary 2.5.

By \cite{superell}, \S 2.6, the conductor exponent $f_2$ has the form
\[
      f_2= \epsilon+\delta,
\]
where
\[
      \epsilon = 2g_Y - \dim H^1_\et(\Zb_k,\QQ_\ell) = 8 - 2 = 6
\]
and $\delta=\delta_V$ is a {\em Swan conductor} of the $\Gamma$-module
$V:=H^1_\et(Y_{\bar{\QQ}},\QQ_\ell)$. Since the graph of $\Yb$ is a tree, the
cospecialization map induces a $\Gamma$-equivariant isomorphism
\[
    V\cong \oplus_{i=0}^3 H^1_\et(\Yb_{i,k},\QQ_\ell).
\]
The Swan conductor of $V$ only depends on the action of $\Gamma_0$. By Lemma
\ref{lem:exa4.2} (iii), the $\Gamma_0$-module $V$ has a direct sum
decomposition
\[
    V = V_0\oplus V_1, 
\]
where
\[
      V_0 = H^1_\et(\Yb_{0,k}), \quad V_1 = \oplus_{i=1}^3 H^1_\et(\Yb_{i,k}).
\]
Moreover, $V_0$ has trivial $\Gamma_0$-action. We conclude that
$\delta=\delta_{V_1}$ is the Swan conductor of the induced $\Gamma_0$-module 
\[
      V_1 = \Ind_{\tilde{\Gamma}_0}^{\Gamma_0} \tilde{V}, \quad
         \tilde{V}:=H^1_\et(\Yb_{1,k},\QQ_\ell),
\]
where the group $\tilde{\Gamma}_0$ is defined in Lemma \ref{lem:exa4.2} (iv).
We have $\delta=\delta_{V_1}=\delta_{\tilde{V}}$ by standard properties of
the Swan conductor. To compute $\delta_{\tilde{V}}$ we may use the formula
\[
   \delta_{\tilde{V}} = \sum_{i=1}^\infty 
     \frac{\abs{\tilde{\Gamma}_i}}{\abs{\tilde{\Gamma}_0}}\cdot
         \dim \tilde{V}/\tilde{V}^{\tilde{\Gamma}_i},
\]
see \cite{superell}, proof of Theorem 2.9. By Lemma \ref{lem:exa4.2} (iv) we
have
\[
     \frac{\abs{\tilde{\Gamma}_i}}{\abs{\tilde{\Gamma}_0}}\cdot
       \dim \tilde{V}/\tilde{V}^{\tilde{\Gamma}_i} = 
       \begin{cases}
           \frac{2}{6}\cdot 2, & i=0,\ldots,15,\\
                            0, & i\geq 16.
       \end{cases} 
\]
We conclude that $\delta=\delta_{\tilde{V}}=10$ and hence 
\[
   f_2=\epsilon+\delta=6+10=16.
\]
\Endproof

It follows that the conductor of the $L$-function is $N= 2^{16} \cdot 317
\approx 10^7$. Using the bound  $M=101248$, we have verified the functional
equation for $L(Y,s)$ and obtained the root number $-1$.

\subsection{}

Finally, we treat a non-hyperelliptic example. Let $Y/\QQ$ be the
superelliptic curve of genus three given by the equation
\begin{equation}  \label{eqn:exa3.1}
  y^3 = f(x) = x^4-x^2 +1.
\end{equation}
The discriminant of $f$ is $144=2^4\cdot 3^2$. We conclude that $Y$ has good
reduction at $p\neq 2,3$. 

The local $L$-factor and the conductor exponent of $Y$ at $p=2$ have been
computed in \cite{superell}, \S~7. The result is
\[
   L_2^{-1} = 1+2T^2, \qquad f_2=8.
\]

The methods of \cite{superell} do not allow the computation of the
semistable reduction of $Y$ at $p=3$, because the exponent of $y$ in
\eqref{eqn:exa3.1} is equal to $p=3$. Again, we have to use the algorithm of
\cite{superp}.

Let $L:=\QQ_3[\zeta_4, \pi]$, where $\zeta_4$ is a primitive $4$th
root of unity and $\pi$ satisfies $\pi^{12}=3$. This is a Galois
extension of $\QQ_3$ whose Galois group is the dihedral group of order
$24$, generated by
\[
\tau(\pi, \zeta_4)=(\zeta_{12}\pi, \zeta_4) \qquad 
\sigma(\pi, \zeta_4)=(\pi, -\zeta_4).
\]
Here $\zeta_{12}:=\zeta_4^3(-1/2+\pi^6\zeta_4/2)\in L$ is a primitive
$12$th root of unity. We also put $\zeta_3=\zeta_{12}^4$.  The residue
field of the unique prime $\p=(\pi)\lhd \OO_L$ is
$\FF_\p=\FF_3[\zeta_4]=\FF_9$.  

\begin{lem}\label{lem:moreexa3}
\begin{enumerate}
\item The curve $Y_L=Y\otimes_{\QQ_3} L$ has semistable reduction at $\p$.
\item Let $\Y$ be the minimal semistable model of $Y$ over $\OO_L$ and
  $\Yb$ the special fiber of $\Y$. Then $\Yb$ consists of $4$ smooth,
  absolutely irreducible components over $\FF_\p$. The normalization
  $\Yb_0$ of the naive model has genus zero. The other three
  components $\Yb_i$ ($i=1,2,3$) have genus $1$, and intersect $\Yb_0$
  in a unique ordinary double point of degree $1$ (Figure
  \ref{fig:Yb_4_3}). The curve $\Yb$ does not have any further
  singularities.
\item The Galois group $\Gamma:=\Gal(L/\QQ_3)$ acts trivially on
  $\Yb_0$. It acts as a cyclic group
  $\langle\psi_{\tau^3}\rangle$ of order $4$ on $\Yb_1$, the quotient
  by this action has genus $0$. 

The components $\Yb_2$ and $\Yb_3$ are conjugate under the action on
$\Yb$ induced by $\sigma.$ The Galois group $\Gamma$ acts on $\Yb_2$
(resp.~$\Yb_3$) as a cyclic group $\langle\psi_\tau\rangle$ of order
$12$. The quotients of $\Yb_2$ (resp.\ $\Yb_3$) both by
$\langle\psi_\tau\rangle$ and by the wild subgroup
$\langle\psi_\tau^4\rangle$ have genus $0$.
\end{enumerate}
\end{lem}

\begin{figure}[h!]\centering 
\begin{tikzpicture}[scale = 0.12]
 \coordinate (A) at (-25,-1);
 \coordinate (D) at (-20,-1.5);
 \coordinate (E) at (-11,4);
 \coordinate (F) at (-13,4);
 \coordinate (G) at (-5,-2);
 \coordinate (K) at (4,4);
 \coordinate (L) at (2,4);
 \coordinate (M) at (11,-2);
 \coordinate (R) at (19,4);
 \coordinate (S) at (17,4);
 \coordinate (T) at (24,-1);
 \coordinate (Z) at (30,-1);

 \coordinate (A1) at (-25,15);
 \coordinate (E1) at (-16,12);
 \coordinate (F1) at (-14,13);
 \coordinate (F2) at (-8,21);
 \coordinate (G1) at (-6,22);
 \coordinate (K1) at (-1,12);
 \coordinate (L1) at (1,13);
 \coordinate (L2) at (7,21);
 \coordinate (M1) at (9,22);
 \coordinate (R1) at (14,12);
 \coordinate (S1) at (16,13);
 \coordinate (S2) at (22,21);
 \coordinate (T1) at (24,22);
 \coordinate (Z1) at (30,15);
% Hilfslinien
%\draw [gray!77] plot coordinates {   (A)(D)(E)(F)(G)(K)(L)(M)(R)(S)(T)(Z)} ;

%singular curve
\draw [black] plot [smooth, tension=1.0] coordinates {  (A)(D)(E)(F)(G)(K)(L)(M)(R)(S)(T)(Z)} ;
\node [right] at (Z) {$\bar Y^{\rm naive}$};

%normalized curve, genus zero part
\draw [black] plot [smooth, tension=1.0] coordinates {  (A1)(Z1)} ;
\node [right] at (Z1) {$\bar Y_0$};

%normalized curve, genus one parts
\draw [black] plot [smooth, tension=1.0] coordinates {  (E1)(F1)(F2)(G1)} ;
\draw [black] plot [smooth, tension=1.0] coordinates {  (K1)(L1)(L2)(M1)} ;
\draw [black] plot [smooth, tension=1.0] coordinates {  (R1)(S1)(S2)(T1)} ;

% delete the loops
 \filldraw [gray!00]      (-12.2,6.7) circle (102pt) ;
 \filldraw [gray!00]      (2.9,7.205) circle (117.5pt) ;
 \filldraw [gray!00]      (18.2,6.56) circle (102pt) ;

\node [below] at (E1) {$\bar Y_1$};
\node [below] at (K1) {$\bar Y_2$};
\node [below] at (R1) {$\bar Y_3$};

\end{tikzpicture}
{\caption{$\Yb\to\Yb^{\rm naive}$}\label{fig:Yb_4_3}}
\end{figure}

\proof Note that $\fb'=\xb(\xb^2+1)$. It follows that the special
fiber $\Yb_k^{\rm naive}$ of the naive model of $Y$ has singularities
in the 3 points with $\xb=0, \pm \zeta_4$. The normalization $\Yb_0$
of $\Yb_k^{\rm naive}$ has genus $0$, since the  map 
\[
(\xb, \yb)\mapsto \xb, \qquad \Yb_k^{\rm naive}\to \PP^1_k
\]
is purely inseparable.

To find the irreducible components $\Yb_i$ for $i=1,2$ we use the coordinates
\[
\begin{split}
x=\pi^9 x_1,\qquad
&y=(\zeta_{3}-1)y_1+1,\\ x=a\pi^{15}x_2-2\zeta_4,\qquad
&y=i\pi^{10}y_2+b\pi^4\left(
-\frac{4}{7}\zeta_4x+1\right),
\end{split}
\]
where $a, b\in \QQ_3[i]$ satisfy $a^2=\zeta_4$ and $b^3=7$.  The
component $\Yb_3$ is obtained by applying the automorphism of $\Yb$
induced by sending $\zeta_4$ to $\zeta_4^3$.  Suitably normalizing $a$
and $b$, we find the following equations:
\[
\begin{split}
  \Yb_1:\qquad &\yb_1^3-\yb_1=\xb_1^2,\\
  \Yb_2:\qquad &\yb_2^3-\yb_2=\xb_2^2.
  \end{split}
  \]
Statements (i) and (ii) follow from this. Statement (iii) follows by a
straightforward verification.
\Endproof

\begin{cor}
  The local $L$-factor and conductor exponent of $Y$ at $p=3$ are
  \[
  L_3^{-1}=1,\qquad f_3=12.
  \]
\end{cor}

\proof The proof is similar to the proof of Corollaries \ref{cor:exa1}
and \ref{cor:exa2}. The statement on the local $L$-factor immediately
follows from Lemma \ref{lem:moreexa3} (iii), since the inertial
reduction $\Zb$ of $Y$ at $p=3$ has genus $0$. (Here we use once more
\cite{superell}, Corollary 2.5.)

We compute the conductor exponent  using
\cite{superell}, \S~2.6. We find that
$f_3=\epsilon+\delta$, where
\[
\epsilon=2g_Y-\dim H^1_\et(\Zb, \QQ_\ell)=6-0=6.
\]

We compute the Swan conductor $\delta$ of the $\Gamma$-module
$V:=H^1_\et(\Zb, \QQ_\ell)$ using \cite{superell}, Theorem 2.9.
We note that the wild subgroup of the
decomposition group of $\p$ satisfies
\[
\tau^4(\pi)-\pi=(\zeta_{12}^4-1)\pi.
\]
Since $v_L(\zeta_{12}^4-1)=6$, we conclude that the filtration of
higher ramification groups is
\[
\Gamma_0=\langle \tau \rangle \supsetneq
\Gamma_1=\cdots=\Gamma_6=\langle\tau^4\rangle\supsetneq
\Gamma_7=\{1\}.
\]

Lemma \ref{lem:moreexa3} (iii) implies that the quotient
$\Yb/\Gamma_i$ has genus one for all $1\leq i\leq 6$. We conclude from
\cite{superell}, Theorem 2.9 that
\[
\delta=\sum_{i=1}^\infty
\frac{|\Gamma_i|}{|\Gamma_0|}\left(2g_Y-2g(\Yb/\Gamma_i)\right)
=\frac{6\cdot 3}{12}(6-2)=6.
\]
It follows  that
\[
f_3=\epsilon+\delta=6+6=12.
\]
\Endproof

Using the bound $M=274994$, we have verified the functional
equation for $L(Y/\QQ)$ and obtained the root number $1$.

\vspace{5ex}\noindent {\small Michel B\"orner, Irene Bouw, Stefan
  Wewers\\ Institut f\"ur Reine Mathematik\\ Universit\"at
  Ulm\\ Helmholtzstr.\ 18\\ 89081 Ulm\\ {\tt
    michel.boerner@uni-ulm.de, irene.bouw@uni-ulm.de,
    stefan.wewers@uni-ulm.de}}

\begin{thebibliography}{10}

\bibitem{ArzdorfDiss}
K.~Arzdorf.
\newblock {\em Semistable reduction of cyclic covers of prime power degree}.
\newblock PhD thesis, Leibniz Universit\"at Hannover, 2012.
\newblock http://edok01.tib.uni-hannover.de/edoks/e01dh12/71609648.pdf.

\bibitem{ArzdorfWewers}
K.~Arzdorf and S.~Wewers.
\newblock Another proof of the semistable reduction theorem.
\newblock Preprint, arXiv:1211.4624, 2012.

\bibitem{AubryPerret}
Y.~Aubry and M.~Perret.
\newblock {Weil} theorems for singular curves.
\newblock In {\em Arithmetic geometry and coding theory (Luminy, 1993)}, pages
  1--7. de Gruyter, Berlin, 1996.

\bibitem{booker}
A.~Booker.
\newblock {\em Numerical tests of modularity}.
\newblock PhD thesis, Princeton University, 2003.

\bibitem{booker05}
A.~Booker.
\newblock Numerical tests of modularity.
\newblock {\em J. Ramanujan Math. Soc.}, 20(4):283--339, 2005.

\bibitem{superell}
I.I. Bouw and S.~Wewers.
\newblock Computing {$L$}-functions and semistable reduction of superelliptic
  curves.
\newblock Preprint, arXiv:1211.4459.

\bibitem{Deligne}
P.~Deligne.
\newblock Les constantes des \'equations fonctionelles des fonctions {$L$}.
\newblock In {\em Modular functions of one variable, II}, number 349 in LNM,
  pages 501--597. Springer-Verlag, 1973.

\bibitem{Dokchitser04}
T.~Dokchitser.
\newblock Computing special values of motivic {$L$-functions}.
\newblock {\em Experimental Mathematics}, 13(2):137--149, 2004.

\bibitem{DokdeJeuZagier}
T.~Dokchitser, R.~de~Jeu, and D.~Zagier.
\newblock Numerical verification of {Beilison's} conjecture for {$K_2$} of
  hyperelliptic curves.
\newblock {\em Compositio Math.}, 142:339--373, 2006.

\bibitem{GaudryHarley2000}
P.~Gaudry and R.~Harley.
\newblock Counting points on hyperelliptic curves over finite fields.
\newblock In W.~Bosma, editor, {\em Algorithmic number theory}, number 1838 in
  Lecture Notes in Computer Science, pages 313--332. Springer-Verlag, 2000.

\bibitem{Kedlaya2001}
K.S. Kedlaya.
\newblock Counting points on hyperelliptic curves using {Monsky--Washnitzer}
  cohomology.
\newblock {\em Journal Ramanujan Math.\ Soc.}, 16:323--338, 2001.

\bibitem{KedlayaSutherland}
K.S. Kedlaya and A.~Sutherland.
\newblock Computing {$L$-series} of hyperelliptic curves.
\newblock In {\em Algorithmic Number Theory}, number 5011 in Lecture Notes in
  Computer Science, pages 312--326. Springer-Verlag, 2008.

\bibitem{liu}
Q.~Liu.
\newblock {\em Algebraic geometry and arithmetic curves}.
\newblock {Oxford University Press}, 2006.

\bibitem{MilneEC}
J.S. Milne.
\newblock {\em \'Etale cohomology}.
\newblock Princeton Univ.\ Press, 1980.

\bibitem{superp}
J.~R\"uth and S.~Wewers.
\newblock Semistable reduction of superelliptic curves of degree $p$.
\newblock in preparation.

\end{thebibliography}
\end{document}